\input amstex
\input epsf
\documentstyle{amsppt}
\magnification=1200
\pageheight{23truecm}
\pretolerance=200
\tolerance=400
\TagsOnRight
\NoRunningHeads
\nologo

\define\dd{\partial}
\define\f{\varphi}
\define\fib{\pmatrix2&1\\ 1&1\endpmatrix}
\define\Int{\operatorname{Int}}
\define\lra{\longleftrightarrow}
\define\N{\Bbb N}
\define\norm#1{\Vert#1\Vert}

\define\RP{\Bbb RP}
\define\sm{\setminus}
\define\Sol{\operatorname{Sol}}
\define\Tor{\operatorname{Tor}}
\define\Tr{\operatorname{Tr}}
\define\Vol{\operatorname{Vol}}
\define\Z{\Bbb Z}

\topmatter
\title Complexity of Torus Bundles Over the Circle With Monodromy $\fib^n$
\endtitle
\author S.~Anisov \endauthor
\address European PostDoctoral Institute \endaddress
\curraddr Dept\. of Mathematics, Utrecht University, P.O\. Box 80{.}010, 
3508~TA Utrecht, the NETHERLANDS	\endcurraddr
\email anisov\@mccme.ru, anisov\@math.uu.nl \endemail
\abstract We find the exact values of complexity for an infinite series of 
3-manifolds. Namely, by calculating hyperbolic volumes, we show that $c(N_n)=
2n$, where $c$ is the complexity of a 3-manifold and $N_n$ is the total space 
of the punctured torus bundle over $S^1$ with monodromy $\fib^n$. We also apply 
a recent result of Matveev and Pervova to show that $c(M_n)\ge2Cn$ with $C
\approx0.598$, where a compact manifold $M_n$ is the total space of the torus 
bundle over $S^1$ with the same monodromy as $N_n$, and discuss an approach to 
the conjecture $c(M_n)=2n+5$ based on the equality $c(N_n)=2n$.
\endabstract
\endtopmatter

\document

The notion of complexity of 3-di\-men\-si\-onal manifolds (see Definition~4 
below) was introduced by S.~Matveev, see~\cite5. Upper bounds for complexity 
can easily be obtained. On the other hand, no lower bounds were known until 
recently, except for only several hundreds of manifolds of small complexity, 
where a full case-by-case analysis can be performed by a computer~\cite{6,~7}; 
thus, neither exact values nor even reasonable lower bounds were known for any 
infinite class of 3-mani\-folds. First meaningful lower bounds of the 
complexity were obtained in~2001, see~\cite8. 

In the present paper, the exact value of complexity is found for an infinite 
series of 3-mani\-folds. This is done in Section~3 for the manifolds $N_n$ that 
are $n$-fold covers of the figure eight knot complement $N_1$; alternatively, 
$N_n$ can be described as the total space of the punctured torus bundle over 
the circle with monodromy $\fib^n$. In Section~2 we apply the aforementioned 
lower bound~\cite8 to the total space $M_n$ of (compact) torus bundle over 
$S^1$ with the same monodromy. Section~1 contains necessary definitions.

\subhead 1. Definitions \endsubhead

In this section, we recall some definitions following~\cite{5,~6}. By $K$ 
denote the 1-di\-men\-si\-onal skeleton of the tetrahedron, which is just the 
clique (that is, the complete graph) with 4 vertices. Note that $K$ is 
homeomorphic to a circle with three radii.

\definition{Definition~1} A compact 2-dimensional polyhedron is called {\it
almost simple\/} if the link of its every point can be embedded in~$K$. An
almost simple polyhedron $P$ is said to be {\it simple\/} if the link of each
point of~$P$ is homeomorphic to either a circle or a circle with a diameter or
the whole graph~$K$. A point of an almost simple polyhedron is {\it
non-singular\/} if its link is homeomorphic to a circle, it is said to be a
{\it triple point\/} if its link is homeomorphic to  a circle with a diameter,
and it is called a {\it vertex\/} if its link is homeomorphic to~$K$. The set
of singular points of a simple polyhedron~$P$ (i.e., the union of the vertices
and the triple lines) is called its {\it singular graph\/} and is denoted
by~$SP$.							\enddefinition

It is easy to see that any compact subpolyhedron of an almost simple polyhedron
is almost simple as well. Neighborhoods of non-singular and triple points of a
simple polyhedron are shown in Fig.~1\,a,\,b; Fig.~1\,c--f represents four
equivalent ways of looking at vertices; in particular, Fig.~1\,e shows the cone 
over the 1-di\-men\-si\-onal skeleton of the tetrahedron.

\midinsert

\epsfxsize=300pt

\centerline{\epsffile{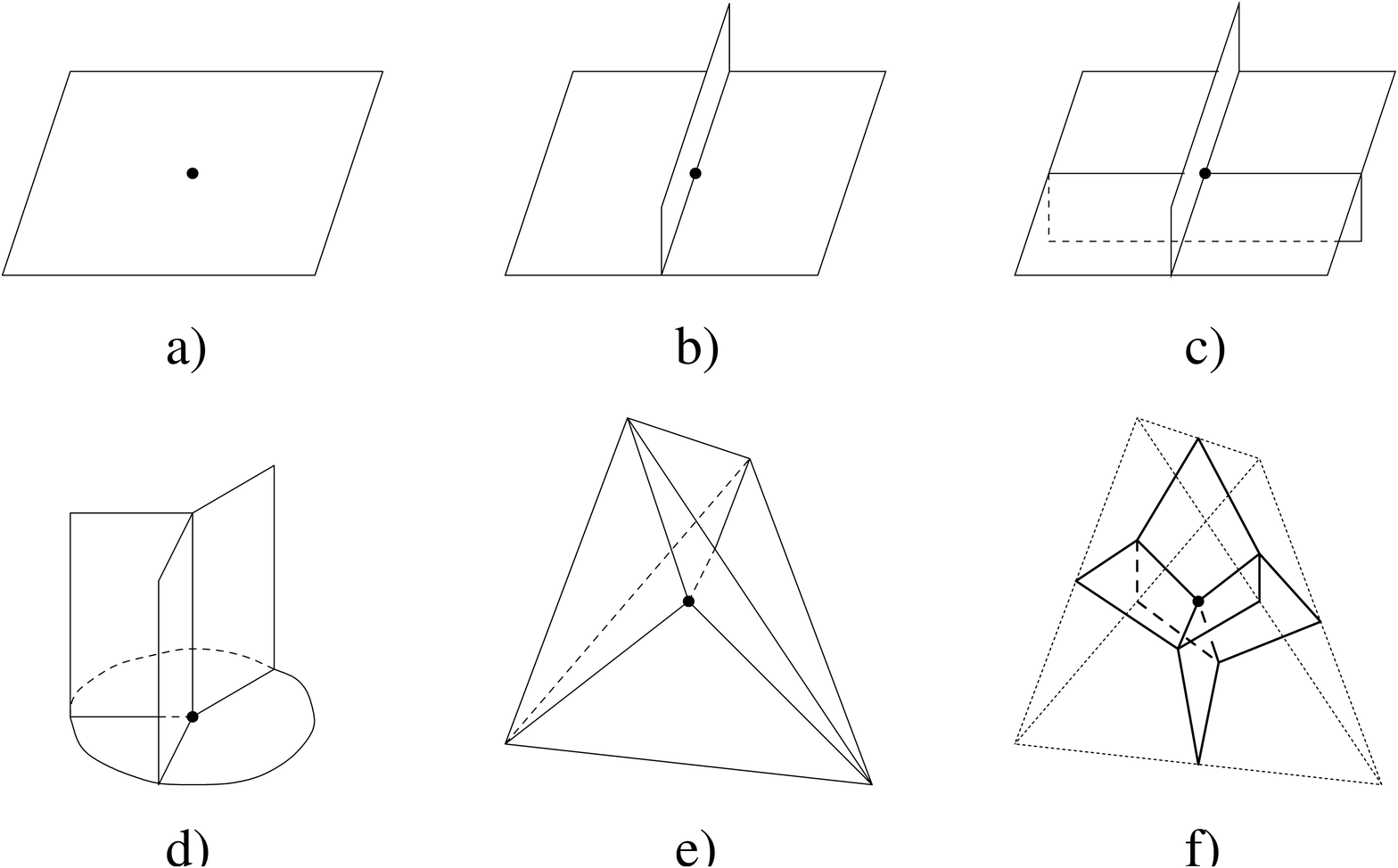}}

\botcaption{Figure 1} Nonsingular (a) and triple (b) points; ways of looking
at vertices (c--f)						\endcaption

\endinsert

\definition{Definition~2} A simple polyhedron~$P$ with at least one vertex is
said to be {\it special\/} if it contains no closed triple lines (without
vertices) and every connected component of $P\sm SP$ is a 2-dimensional cell.
								\enddefinition

\definition{Definition~3} A polyhedron $P\subset\Int M$ is called a {\it
spine\/} of a compact 3-di\-men\-sional manifold~$M$ if $M\sm P$ is
homeomorphic to $\dd M\times(0,1]$ (if $\dd M\ne0$) or to an open 3-cell (if 
$\dd M=0$). In other words, $P$ is a spine of~$M$ if a manifold $M$ with 
boundary (or a closed manifold~$M$ punctured at one point) can be collapsed 
onto~$P$. A spine $P$ of a 3-manifold $M$ is said to be {\it almost simple}, 
{\it simple}, or {\it special\/} if it is an almost simple, simple, or special 
polyhedron, respectively.					\enddefinition

Given a special spine $P$ of a compact manifold $M^3$, one can construct a dual 
singular triangulation of $M^3$ with one vertex (lying in the middle of the 
3-cell $M\sm P$), see Fig.~1\,f; if $M$ is a manifold with connected boundary, 
the same construction gives a triangulation of the one-point compactification 
of $M\sm\dd M$. In both cases, there is a one-to-one correspondence between 
vertices of~$P$ and tetrahedra of the triangulation.

\definition{Definition~4} The {\it complexity\/} $c(M)$ of a compact
3-manifold~$M$ is the minimal possible number of vertices of an almost simple
spine of~$M$. An almost simple spine with the smallest possible number of
vertices is said to be a {\it minimal\/} spine of~$M$.		\enddefinition

\proclaim{Theorem~1~\cite5} Let $M$ be an orientable irreducible $3$-manifold 
with incompressible \rom(or empty\rom) boundary and without essential annuli. 
If $c(M)>0$ \rom(that is, if $M$ is different from \rom(possibly punctured\rom) 
$S^3$, $\RP^3$, and $L_{3,1}$\rom), then any minimal almost simple spine of~$M$ 
is special. \endproclaim

Thus, if $M$ is as in Theorem~1, then $c(M)$ is equal to the minimal number of 
tetrahedra in a singular triangulation of~$M$. By the way, this implies that 
$c(M)\ge\norm M$, where $\norm M$ stands for the Gromov norm of~$M$ 
(see~\cite4), whenever $M$ is a compact 3-manifold that satisfies the 
assumptions of Theorem~1. We do not know any manifold~$M$ such that 
$c(M)=\norm M$.

\subhead 2. Torus bundles over $S^1$: a lower bound \endsubhead

\proclaim{Theorem~2~\cite8} Let $M$ be a compact irreducible orientable 
\rom3-manifold non-homeo\-morphic to $S^3$, $\RP^3$, and~$L_{3,1}$. Then 
$c(M)\ge2\log_5\left|\Tor(H_1(M,\Z))\right|+\beta_1(M,\Z)-1$.	\endproclaim

\proclaim{Corollary~1~\cite8} For lens spaces, we have $c(L_{p,q})\ge2
\log_5p-1$.							\endproclaim

In particular, $c(L_n)\ge2\log_5\f_n-1$, where $L_n$ stands for $L_{\f_n,
\f_{n-1}}$ and $\f_n$ denotes the $n$\snug th Fibonacci number. Since $\f_n=(
((\sqrt5+1)/2)^n+((-\sqrt5+1)/2)^n)/\sqrt5$, we have $c(L_n)\ge C_nn-2$ with 
$C_n=\frac2n\log_5(\sqrt5\f_n)$ tending 
to $C=2\log_5((\sqrt5+1)/2)\approx0.598$ as $n\to\infty$, which is a fairly 
good estimate, since $c(L_n)\le n-4$ whenever $n\ge4$, see~\cite6.

Theorem~2 can be successfully applied to some other 3-manifolds. Let us denote
by~$M_n$ the total space of the $T^2$-bundle over~$S^1$ with monodromy
$A^n$, where $A=\fib$. A short calculation shows that $\left|\Tor(H_1(M_n,\Z))
\right|=\pm\det(A^n-I)=\pm(\det A^n-\Tr A^n+1)$; since $A^n=\pmatrix\f_{2n+1}&
\f_{2n}\\ \f_{2n}&\f_{2n-1}\endpmatrix$ and $\det A^n=1$, we have $\left|\Tor
(H_1(M_n,\Z))\right|=\f_{2n+1}+\f_{2n-1}-2$. Taking into account that $\beta_1
(M_n,\Z)=1$, we get the following estimate in a similar way.

\proclaim{Corollary~2} $c(M_n)\ge2C_nn$, where $C_n=\frac1n\log_5(\f_{2n+1}+
\f_{2n-1}-2)$.							\endproclaim

Note that $C_n\to C=\log_5((\sqrt5+1)/2)^2\approx0.598$ as $n\to\infty$, which
is as good as in the previous example since $c(M_n)\le2n+5$, see~\cite1.
Combining this with the inequalities $c(M_n)\ge7$ and $C_n>0.597$ whenever $n
\ge6$ (none of the $M_n$ is contained in the list of 3-manifolds up to 
complexity~6, see~\cite6; in fact, all compact 3-manifolds up to complexity~6 
are elliptic except for the flat manifolds, which all have complexity~6, while 
all the manifolds $M_n$ are $\Sol$-manifolds), we get $c(M_n)>1.19n$ for 
all~$n\ge1$. We believe that $c(M_n)=2n+5$ (see~\cite1) and $c(L_n)=n-
4$~(see~\cite{5,~6}). Note that $\norm{M_n}=0$, because there is an obvious 
action of~$S^1$ on~$M_n$.

\subhead 3. Punctured torus bundles: exact values \endsubhead

This is the main section of the paper. Here we find the exact values of $c(N_n
)$ for an infinite series of 3-manifolds; to the best of our knowledge, this is 
the first result of this kind. The manifolds $N_n$, $n\in\N$, are the total 
spaces of the punctured torus bundles over~$S^1$ with monodromy~$A^n$, where 
$A=\fib$.

The manifold $N_1$ has been studied extensively. It is well known to be 
hyperbolic and homeomorphic to the figure eight knot complement, 
see~\cite{9,~Chap.~4}. A special spine $P$ of~$N_1$ with two vertices, four 
edges, and two hexagonal 2-cells is represented by Fig.~2\,a; the picture shows 
the boundary of the neighborhood in~$P$ of the singular graph~$SP$ (which 
consists of two vertices and four triple lines; note that the triple lines 
themselves are not drawn). Since the spine~$P$ is special (and thus both 
components of $P\sm SP$ are disks), the picture contains enough information to 
reconstruct~$P$. This spine coincides with the spine $P_0$ constructed
in~\cite{1,~\S2.3}. Figures~2\,a and~2\,c are taken from~\cite{2,~\S4}; many 
other pictures related to the manifold~$N_1$ can be found in~\cite{3,~Chap.~8}.

\midinsert

\epsfxsize=360pt

\centerline{\epsffile{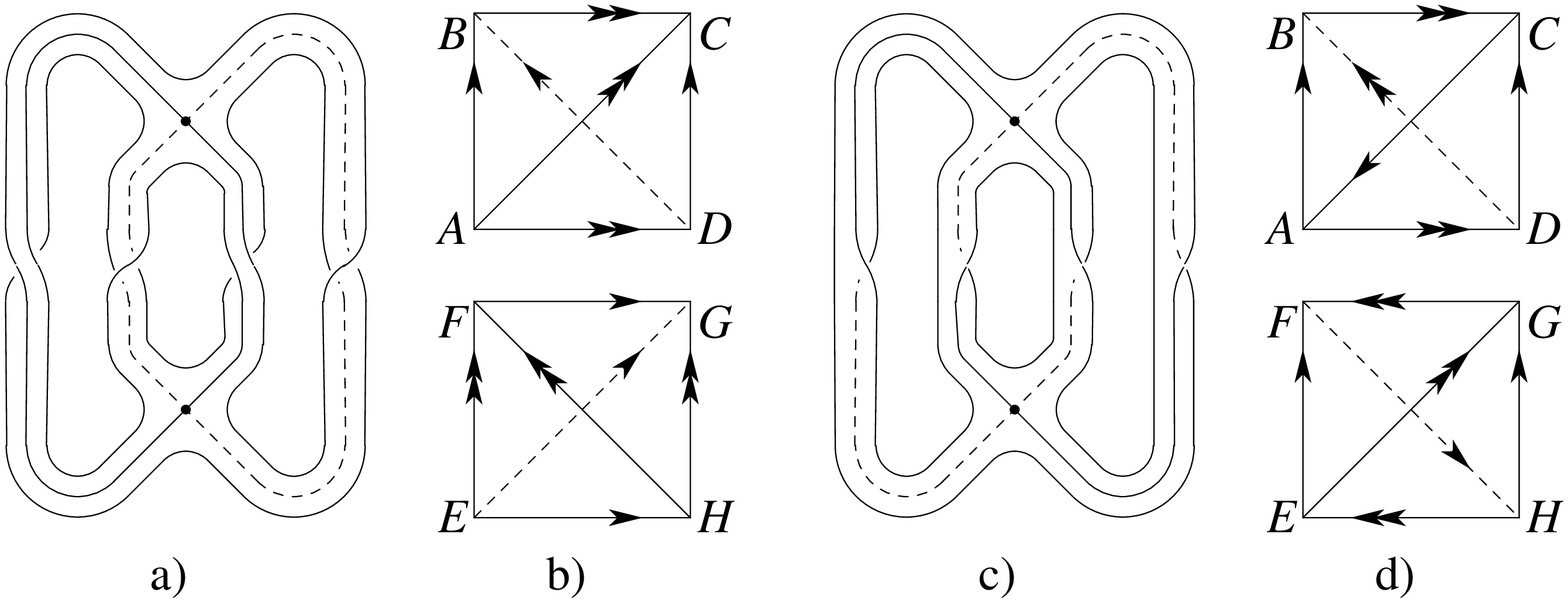}}

\botcaption{Figure 2} Noncompact hyperbolic 3-manifolds of complexity~2
\endcaption

\endinsert

\proclaim{Theorem~3} The equality $c(N_n)=2n$ holds.		\endproclaim

\demo{Proof} Any special spine defines its dual decomposition of the 
manifold into tetrahedra. For $N_1$, this decomposition is shown on Fig.~2\,b. 
It consists of two tetrahedra (each contains one vertex of the spine), glued 
together so that there are two edges (dual to the hexagons of the spine); each 
of them is incident to six dihedral angles of the tetrahedra; it is described 
in detail in~\cite{9,~Chap.~1}. The gluing pattern can be reconstructed from 
Fig.~2\,a; it is $ABC\lra EHF$, $BAD\lra GEF$, $CDA\lra GFH$, $DCB\lra EGH$. 
Thus, the edges marked by single arrows are glued together, those marked by 
double arrows are glued together, too, and the direction of the arrows is 
respected. 

A complete hyperbolic structure on $N_1$ comes from that on two tetrahedra 
considered as regular ideal tetrahedra in~$H^3$. All their dihedral angles are
equal to~$\pi/3$~\cite9, so the sum of the dihedral angles incident to an edge 
equals $2\pi$ for both edges, which means that the hyperbolic structure 
described above is well defined. Among the ideal tetrahedra in~$H^3$, the 
regular one has the maximal volume $V\approx1.0149$, see, 
e.g.,~\cite{9,~Chap.~7}. Therefore, the manifold $N_1$ admits a hyperbolic 
structure of volume~$2V$.

Since there is an $n$-fold covering $p\:N_n\to N_1$, the polyhedron $p^{-1}(P)$ 
is a special spine of $N_n$ with $2n$ vertices, so we have $c(N_n)\le 2n$; 
again, that spine coincides with the one constructed in~\cite{1,~\S2.3}. For 
the same reason, the manifold $N_n$ admits a complete hyperbolic structure of 
volume~$2nV$. Now we have to prove the inequality $c(N_n)\ge 2n$. 

The manifolds $N_n$ satisfy the hypotheses of Theorem~1. Thus, their minimal 
spines are special. So, if a minimal spine of~$N_n$ contains $k$ vertices, then 
there is a (singular) triangulation of $N_n$ formed by $k$ tetrahedra. 
Straightening them, we get a triangulation of a fundamental domain for $\pi_1(
N_n)$ in~$H^3$, which has volume~$2nV$, into $k$~ideal tetrahedra (which may 
overlap). Since the volume of any ideal tetrahedron in~$H^3$ does not 
exceed~$V$, we get~$k\ge2n$.					\qed\enddemo

\remark{Remarks} 1.~In fact, we have shown that $c(M)\ge\left\lceil\frac{\Vol
(M)}V\right\rceil$ for any hyperbolic manifold $M^3$, either compact or 
noncompact, orientable or not. We know no examples of compact hyperbolic 
3-mani\-folds for which this estimate is sharp. On the other hand, there exist 
compact orientable hyperbolic 3-mani\-folds of volume $0.94...$ and $0.98...$, 
while the complexity of any compact orientable hyperbolic 3-mani\-fold is at 
least~9, see~\cite2. Moreover, there exist infinitely many compact hyperbolic 
3-mani\-folds such that their volume is less than~$2V$~\cite{2,~9}; their 
list contains manifolds of arbitrary large complexity, because there are only 
finitely many irreducible 3-manifolds of complexity bounded by any integer~$N$, 
see~\cite5.

2.~There exists one more noncompact orientable hyperbolic 3-manifold of 
volume~$2V$ and complexity~2, see~\cite2. Its minimal special spine and 
corresponding triangulation are shown on Fig.~2\,c,\,d. The gluing pattern is 
$ABC\lra FHE$, $BAD\lra FEG$, $CDA\lra HFG$, $DCB\lra HGE$. The complexity of 
any $n$-fold covering space of this manifold is again equal to~$2n$. The proof 
of this statement repeats that of Theorem~3.

3.~Let us return to the manifolds~$M_n$ considered in Section~2. Consider a 
minimal triangulation of~$M_n$ (dual to its minimal spine~$P$, which is special 
by Theorem~1). Since there is only one vertex (dual to the 3-cell $M\sm P$), 
all the edges of the triangulation are loops. They generate the group~$\pi_1
(M)$. Therefore, at least one of them has a nonzero image under the projection 
$p_*\:\pi_1(M_n)\to\pi_1(S^1)$. Let us suppose for a moment that there is an 
edge~$e$ that is isotopic to the section of the fibration $p\:M_n\to S^1$. Let 
$\sigma$ be the 2-component of~$P$ dual to~$e$. Put~$P'=P\sm\sigma$. Then $P'$ 
is an almost simple spine of the manifold~$M_n\sm e=N_n$. By Theorem~3, 
$P'$~contains at least $2n$ vertices. Consequently, $P=P'\cup\sigma$ has at 
least $2n+2$ vertices, which is close to the conjectured value $c(M_n)=2n+5$, 
see~\cite1 (indeed, if adding $\sigma$ to $P'$ does not increase the number of 
vertices, then $\dd\sigma$ is a closed triple line and the spine~$P$ is not 
minimal by virtue of Theorem~1; if all vertices of~$P$ belonging to~$\dd\sigma$ 
are different but their number is less than~$4$, then a simplification 
move~\cite{5--7} can be applied, and $P$ is not a minimal spine; finally, one 
can show that the case where~$\dd\sigma$ passes through some vertex of~$P$ more 
than once but does not pass through any other vertex is impossible). However, 
it remains unclear why such an edge~$e$ should exist in a triangulation dual to 
arbitrary minimal spine of~$M_n$. \endremark

\definition{Acknowledgements} The main idea of this paper has appeared during 
my visit to Institut Joseph Fourier (Grenoble, France). The paper has been 
finished at Utrecht University, the Netherlands. The author thanks both 
institutes for their kind hospitality. The author has pleasure to thank 
S.~Matveev for useful discussions.
\enddefinition

\enddocument
\end